\documentclass{amsart}

\usepackage{fancyhdr}
\usepackage{lastpage}
\usepackage{stmaryrd,yhmath}

\newcommand{\bdis}{\begin{displaymath}}
\newcommand{\edis}{\end{displaymath}}
\newcommand{\be}{\begin{equation}}
\newcommand{\ee}{\end{equation}}
\newcommand{\mbb}{\mathbb}
\newcommand{\mcal}{\mathcal}

\newcommand{\vp}{\varphi}

\newcommand{\vth}{\vartheta}

\newcommand{\mT}{\mathring{T}}

\newcommand{\zf}{\zeta\left(\frac{1}{2}+it\right)}

\DeclareMathOperator{\sn}{sn}
\DeclareMathOperator{\cn}{cn}
\DeclareMathOperator{\dn}{dn}

\newtheorem{theorem}{Theorem}

\theoremstyle{definition}

\newtheorem{example}[theorem]{Example}

\theoremstyle{remark}
\newtheorem{remark}[]{Remark}

\newtheorem*{mydef1}{{\bf Theorem}}

\newtheorem*{mydef4}{{\bf Corollary}}

\numberwithin{equation}{section}

\begin{document}

\title{Jacob's ladders and $|\zeta|^{-2}$-representation of some functionals generated by a new class of transcendental integrals}

\author{Jan Moser}

\address{Department of Mathematical Analysis and Numerical Mathematics, Comenius University, Mlynska Dolina M105, 842 48 Bratislava, SLOVAKIA}

\email{jan.mozer@fmph.uniba.sk}


\begin{abstract}
In this paper we introduce a new infinite set of transcendental integrals. Each of them is expressed by corresponding value of the function
$|\zf|^{-2}$. Such a property is another argument about universality of the Riemann zeta-function $\zf$ in the field of pure
mathematics.
\end{abstract}
\maketitle

\section{Introduction}

\subsection{}

In our paper \cite{7} we have postulated the radius of the Universe $R(t)$ as
\be \label{1.1}
R(t)=h\left(\left|\zf\right|\right),\ t\to\infty,
\ee
(see \cite{7}, (3.5), (4.1)). By making use this together with the Riemann hypothesis we have constructed:
\begin{itemize}
\item[(a)] a new infinite set of mathematical universes (of the Einstein's type),
\item[(b)] an infinite subset of microscopic universes (on the Planck's scale) of the Einstein's type in the early period (just after inflation)
of the evolution of the Universe.
\end{itemize}

\begin{remark}
Above mentioned results can be understood as an argument for the universality (comp. \cite{2}, p. 135) of the modulus of Riemann zeta-function
on the critical line.
\end{remark}

\subsection{}

In this paper we will study the property of universality of the function
\bdis
\left|\zf\right|
\edis
in another direction. Namely, we will study a new class of transcendental integrals of the following type
\bdis
\begin{split}
& \int_a^b |H[\vp_1(t);\tau]|^\alpha |H'_{\vp_1}[\vp_1(t);\tau]|{\rm d}t, \\
& a<\tau<b,\quad \alpha>0 ,
\end{split}
\edis
where
\be \label{1.2}
\vp_1(t)\sim t-(1-c)\pi(t),\quad \pi(t)\sim \frac{t}{\ln t},\ t\to\infty,
\ee
and $\pi(t)$ is the prime-counting function and $c$ is the Euler constant. We will prove that there are numbers
\bdis
\bar{a}, \bar{\tau},\bar{b},\bar{\alpha},t_H,
\edis
and a function
\bdis
h_1(u)
\edis
such that
\be \label{1.3}
\begin{split}
& \int_{\bar{a}}^{\bar{b}}|H[\vp_1(t);\tau]|^{\bar{\alpha}}|H'_{\vp_1}(t);\tau|{\rm d}t=\\
& =h_1\left(\left|\zeta\left(\frac 12+it_H\right)\right|\right),\ t_H=t(H;\bar{\alpha}),\ \bar{a}\to\infty
\end{split}
\ee
(comp. with (1.1)), for every fixed Jacob's ladder $\vp_1(t)$.

\subsection{}

The function $\vp_1(t)$ (see (1.1) -- (1.3)) that we call Jacob's ladder (see \cite{5}) according to the Jacob's dream in Chumash,
Bereishis, 28:12 has the following properties:
\begin{itemize}
\item[(a)]
\bdis
\vp_1(t)=\frac 12\vp(t),
\edis
\item[(b)] the function $\vp(t)$ is the solution of nonlinear integral equation (see \cite{5}, \cite{6})
\bdis
\int_0^{\mu[x(T)]}Z^2(t)e^{-\frac{2}{x(T)}t}{\rm d}t=\int_0^TZ^2(t){\rm d}t,
\edis
where
\be \label{1.4}
\begin{split}
& Z(t)=e^{i\vth(t)}\zf , \\
& \vth(t)=-\frac t2\ln\pi+\text{Im}\ln\Gamma\left(\frac 14+i\frac t2\right),
\end{split}
\ee
and each advisable function $\mu(y)$ generates a solution
\bdis
y=\vp_\mu(T)=\vp(T);\ \mu(y)\geq 7y\ln y.
\edis
\end{itemize}

\begin{remark}
The main reason for the introduction of the Jacob's ladders in the paper \cite{5} lies in the proof of the following theorem: \\

the Hardy-Littlewood integral (1918)
\bdis
\int_0^T\left|\zf\right|^2{\rm d}t
\edis
has -- in addition to the Hardy-Littlewood (and other similar) expression possessing an unbounded error term as $T\to \infty$ -- the following
infinite set of almost exact expressions
\be \label{1.5}
\begin{split}
& \int_0^T\left|\zf\right|^2{\rm d}t=\vp_1(T)\ln\vp_1(T)+(c-\ln 2\pi)\vp_1(T)+c_0+\\
& +\mcal{O}\left(\frac{\ln T}{T}\right),\quad T\to\infty,
\end{split}
\ee
where $c_0$ is the constant from the Titchmarsh-Kober-Atkinson formula (see \cite{8}, p. 141).
\end{remark}

\begin{remark}
Simultaneously with (\ref{1.5}) we have proved that the following transcendental equation
\bdis
\int_0^T\left|\zf\right|^2{\rm d}t=V(T)\ln V(T)+(c-\ln 2\pi)V(T)+c_0
\edis
has an infinite set of asymptotic solutions
\bdis
V(T)=\vp_1(T),\ T\to\infty.
\edis
\end{remark}

\section{The result}

The following theorem holds true in the direction (1.3).

\begin{mydef1}
Let
\begin{itemize}
\item[(a)]
\bdis
G(t)\in C_1[T_0,\infty),\ T_0>0,
\edis
\item[(b)] the symbols
\bdis
\{\gamma\}, \{ t_0\};\ t_0\not=\gamma;\ \gamma=\gamma[G], t_0=t_0[G]
\edis
denote the sequence of the roots of equations
\bdis
G(t)=0,\ G'(t)=0
\edis
respectively,
\item[(c)] the points of the sequences $\{\gamma\},\{t_0\}$ are separated each from other, i. e.
\bdis
\gamma'<t_0<\gamma'';\ \gamma''-\gamma'\in \left(\left. 0,\frac{\gamma'}{\ln\gamma'}\right.\right],
\edis
where $\gamma',\gamma''$ are neighbouring points of the sequence $\{\gamma\}$, (i. e. $G(t_0)$ is the local extreme of the function $G(t), t\in [\gamma',\gamma'']$),
\item[(d)]
\bdis
\vp_1\left\{[\mathring{\gamma}',\mathring{\gamma}'']\right\}=[\gamma',\gamma''],
\edis
\item[(e)]
\bdis
\begin{split}
& H[\vp_1(t);t_0]=\frac{G[\vp_1(t)]}{G[\vp_1(\mathring{t}_0)]}=\frac{G[\vp_1(t)]}{G(t_0)}; \ t_0=\vp_1(\mathring{t}_0).
\end{split}
\edis
\end{itemize}
Then there is the function
\be \label{2.1}
\omega(x)=\omega(x;H)=1+\mcal{O}\left(\frac{\ln\ln x}{\ln x}\right),\ x\to\infty
\ee
and the point
\bdis
t_H=t(H;\omega)\in \left(\mathring{\gamma}',\mathring{\gamma}''\right)
\edis
such that
\be \label{2.2}
\begin{split}
& \int_{\mathring{\gamma}'}^{\mathring{\gamma}''} |H[\vp_1(t);t_0]|^{2\omega(\gamma')\ln\gamma'-1}|H'_{\vp_1}[\vp_1(t);t_0]|{\rm d}t= \\
& = \frac{1}{\left|\zeta\left(\frac 12+it_H\right)\right|^2},\ \gamma'\to\infty,
\end{split}
\ee
where, of course,
\bdis
2\omega(\gamma')\ln\gamma'-1\to\infty,\ \gamma'\to\infty.
\edis
\end{mydef1}

\begin{remark}
The formula (\ref{2.2}):
\begin{itemize}
\item[(a)] introduces a new infinite set of transcendental integrals each of them is expressed by the corresponding value of the function
\bdis
\left|\zf\right|^{-2},
\edis
i. e. $h_1(u)=u^{-2}$, comp. (\ref{1.3}),
\item[(b)] supports the argument about the universality of the Riemann zeta-function $|\zf|$ in the field of pure mathematics (comp. Remark 1),
\item[(c)] is not reachable by recent methods in the theory of the Riemann zeta-function.
\end{itemize}
\end{remark}

\section{Examples concerning Jacobi elliptical functions}

\subsection{}

Let us consider the continuum set of functions
\bdis
\sn x=\sn(x;k),\ x>0,\ k^2\in (0,1),
\edis
where (comp. \cite{9}, pp. 35-40)
\bdis
K=K(k)=\int_0^1 \frac{{\rm d}t}{\sqrt{(1-t^2)(1-k^2t^2)}}.
\edis
Then we have
\bdis
\begin{split}
& \gamma=\gamma[\sn];\ \gamma_l=2lK,\ l\in\mbb{N}_0,\ \gamma''-\gamma'=2K=O(1);\ \gamma_l\to\infty, \\
& t_0=t_0[\sn];\ t_0(l)=(2l+1)K, \\
& \sn'x=\cn x\dn x;\ \dn x>0, \\
& |\sn t_0|=|\sn(2l+1)K|=1 \ \Rightarrow \ |G(t_0)|=1 \ \Rightarrow \\
& \Rightarrow\ |H[\vp_1(t);t_0]|=|G[\vp_1(t)]|=|\sn[\vp_1(t)]|.
\end{split}
\edis

Consequently we obtain (see (\ref{2.2})) the following

\begin{example}
There is
\bdis
t_s\in \left( \mathring{\gamma}'[\sn],\mathring{\gamma}''[\sn]\right)
\edis
such that
\bdis
\begin{split}
& \int_{\mathring{\gamma}'}^{\mathring{\gamma}''}|\sn[\vp_1(t)]|^{2\omega(\gamma')\ln\gamma'-1}|\cn[\vp_1(t)]|\dn[\vp_1(t)]{\rm d}t= \\
& = \frac{1}{\left|\zeta\left(\frac 12+it_s\right)\right|^2},\ \gamma'\to\infty.
\end{split}
\edis
\end{example}

\subsection{}

Next, for continuum set of functions
\bdis
\cn x=\cn (x;k),\ x>0
\edis
we have by the similar way that
\bdis
\begin{split}
& \gamma=\gamma[\cn];\ \gamma_l=(2l+1)K,\ l\in\mbb{N}_0, \\
& t_0=t_0[\cn];\ t_0(l)=(2l+2)K, \\
& \cn'x=-\sn x\dn x, \\
& |\cn t_0|=|\cn(2l+2)K|=1 \ \Rightarrow \ |G(t_0)|=1 \ \Rightarrow \\
& \Rightarrow\ |H[\vp_1(t);t_0]|=|G[\vp_1(t)]|=|\cn[\vp_1(t)]|.
\end{split}
\edis
Consequently, we have the following

\begin{example}
There is
\bdis
t_c\in \left( \mathring{\gamma}'[\cn],\mathring{\gamma}''[\cn]\right)
\edis
such that
\bdis
\begin{split}
& \int_{\mathring{\gamma}'}^{\mathring{\gamma}''}|\cn[\vp_1(t)]|^{2\omega(\gamma')\ln\gamma'-1}|\sn[\vp_1(t)]|\dn[\vp_1(t)]{\rm d}t= \\
& = \frac{1}{\left|\zeta\left(\frac 12+it_c\right)\right|^2},\ \gamma'\to\infty.
\end{split}
\edis
\end{example}

\section{Example concerning the Bessel's functions}

Let
\bdis
\gamma=\gamma\left[\frac{J_\nu(x)}{x^\nu}\right],\ t_0=t_0\left[\frac{J_\nu(x)}{x^\nu}\right].
\edis
For our purpose it is sufficient to consider the case
\bdis
\nu>-1,\quad x\to\infty.
\edis
We have from the theory of the Bessel's functions that
\bdis
\gamma'<t<\gamma'',\ \gamma''-\gamma'\sim \pi,\ \gamma'\to\infty,
\edis
(comp. \cite{11}, p. 361; \cite{10}, p. 91) where, of course,
\bdis
\frac{J_{\nu+1}(x)}{x^\nu}=-\frac{{\rm d}}{{\rm d}x}\left[\frac{J_\nu(x)}{x^\nu}\right].
\edis

Next we have
\bdis
\begin{split}
& G[\vp_1(t)]=\frac{J_\nu[\vp_1(t)]}{[\vp_1(t)]^\nu},\ G(t_0)=\frac{J_\nu(t_0)}{t_0^\nu};\ \vp_1(\mathring{t}_0)=t_0, \\
& H[\vp_1(t);t_0]=\frac{J_\nu[\vp_1(t)]}{[\vp_1(t)]^\nu}\frac{t_0^\nu}{J_\nu(t_0)}.
\end{split}
\edis
Consequently, we obtain (see (\ref{2.2})) the following

\begin{example}
There is
\bdis
t_J\in \left(\mathring{\gamma}'\left[\frac{J_\nu(x)}{x^\nu}\right],\mathring{\gamma}''\left[\frac{J_\nu(x)}{x^\nu}\right]\right)
\edis
such that
\bdis
\begin{split}
& \int_{\mathring{\gamma}'}^{\mathring{\gamma}''}
\left|\frac{J_\nu[\vp_1(t)]}{[\vp_1(t)]^\nu}\frac{t_0^\nu}{J_\nu(t_0)}\right|^{2\omega(\gamma')\ln\gamma'-1}
\left| \frac{J_{\nu+1}[\vp_1(t)]}{[\vp_1(t)]^\nu}\frac{t_0^\nu}{J_\nu(t_0)}\right|{\rm d}t= \\
& = \frac{1}{\left|\zeta\left(\frac 12+it_J\right)\right|^2},\ \gamma'\to\infty.
\end{split}
\edis
\end{example}

\section{Example concerning the Riemann zeta-function}

In the case of the function $Z(t)$ (see (\ref{1.4})) we have, on the Riemann hypothesis, that
\bdis
\gamma'<t_0<\gamma'',\ \gamma'\to\infty,
\edis
where
\bdis
\gamma=\gamma[Z],\ t_0=t_0[Z],
\edis
i. e. the points of the sequences $\{\gamma\},\{t_0\}$ are separated each from other (see \cite{4}, Corollary 3). Next, on the Riemann hypothesis also the Littlewood's
estimate
\bdis
\gamma''-\gamma'<\frac{A}{\ln\ln\gamma'},\ \gamma'\to\infty
\edis
holds true (see \cite{3}). Since
\bdis
H[\vp_1(t);t_0]=\frac{Z[\vp_1(t)]}{Z(t_0)},
\edis
then we obtain the following

\begin{example}
On the Riemann hypothesis there is
\bdis
t_Z\in \left( \mathring{\gamma}'[Z],\mathring{\gamma}''[Z]\right)
\edis
such that
\be \label{5.1}
\begin{split}
& \int_{\mathring{\gamma}'}^{\mathring{\gamma}''}
\left|\frac{Z[\vp_1(t)]}{Z(t_0)}\right|^{2\omega(\gamma')\ln\gamma'-1}
\left| \frac{Z'_{\vp_1}[\vp_1(t)]}{Z(t_0)}\right|{\rm d}t= \\
& = \frac{1}{\left|\zeta\left(\frac 12+it_Z\right)\right|^2},\ \gamma'\to\infty.
\end{split}
\ee
\end{example}

\begin{remark}
Consequently, we obtain by (\ref{5.1}) by the first formula in (\ref{1.4}) the following sufficiently complicated integral
\bdis
\begin{split}
& \int_{\mathring{\gamma}'}^{\mathring{\gamma}''}
\left|\frac{\zeta\left[\frac 12+i\vp_1(t)\right]}{\zeta\left(\frac 12+it_0\right)}\right|^{2\omega(\gamma')\ln\gamma'-1}\times \\
& \times \left|\vth'_{\vp_1}[\vp_1(t)]\frac{\zeta\left[\frac 12+i\vp_1(t)\right]}{\zeta\left(\frac 12+it_0\right)}+
\frac{\zeta'_{\vp_1}\left[\frac 12+i\vp_1(t)\right]}{\zeta\left(\frac 12+it_0\right)}\right|{\rm d}t= \\
& = \frac{1}{\left|\zeta\left(\frac 12+it_Z\right)\right|^2},\ \gamma'\to\infty.
\end{split}
\edis
\end{remark}

\section{Topological deformations of some element and $|\zeta|^{-2}$-representation of a corresponding
functionals connected with (\ref{2.2})}

Let some element
\bdis
\bar{G}(t)=C_1 [\gamma',\gamma'']
\edis
fulfil the assumptions (b) and (c) of Theorem. Let
\bdis
\{ G_T(t)\}
\edis
denote the continuum set of all topological deformations $G_T(t)$ of the graph of the function $\bar{G}(t)$ such that every
\bdis
G_T(t),\ t\in [\gamma',\gamma'']
\edis
fulfils (b) and (c) of Theorem, and
\bdis
G_T(t)\in C_1[\gamma',\gamma''].
\edis
Let, finally,
\be \label{6.1}
\begin{split}
& H_T[\vp_1(t);t_0]\in \mbb{D}\{ [\mathring{\gamma}',\mathring{\gamma}'']\} \ \Leftrightarrow \\
& \Leftrightarrow\ H_T[\vp_1(t);t_0]=\frac{G_T[\vp_1(t)]}{G_T[\vp_1(\mathring{t}_0)]},\ t\in[\mathring{\gamma}',\mathring{\gamma}''],
\end{split}
\ee
where
\bdis
\vp_1\{[\mathring{\gamma}',\mathring{\gamma}'']\}=[\gamma',\gamma''].
\edis
Then we have that by the integral
\be \label{6.2}
\begin{split}
& \mcal{F}[H_T]=\int_{\mathring{\gamma}'}^{\mathring{\gamma}''}
|H_T[\vp_1(t);t_0]|^{2\omega(\gamma')\ln\gamma'-1}\times \\
& \times |(H_T)'_{\vp_1}[\vp_1(t);t_0]|{\rm d}t,\ H_T\in\mbb{D}
\end{split}
\ee
is defined the functional on $\mbb{D}$.

\begin{mydef4}
The formula (see (\ref{2.2}), (\ref{6.2}))
\bdis
\begin{split}
& \mcal{F}[H_T]=\frac{1}{\left|\zeta\left(\frac 12+it_H\right)\right|^2},\ H_T\in\mbb{D}, \\
& t_{H_T}\in (\mathring{\gamma}',\mathring{\gamma}''),\ \gamma'\to\infty
\end{split}
\edis
expresses the representation of the functional $\mcal{F}[H_T]$ by the corresponding set of values of the function $|\zf|^{-2}$.
\end{mydef4}

\section{Proof of Theorem}

Let us remind that we have proved the following theorem (see \cite{6}, (9.7)): for every Lebesgue integrable function
\bdis
\begin{split}
& f(x),\ x\in [T,T+U],\quad f(x)\geq 0\ (\leq 0),\ x\in [T,T+U]; \\
& \vp_1\{[\mT,\widering{T+U}]\}=[T,T+U]
\end{split}
\edis
we have
\be \label{7.1}
\begin{split}
& \int_{\mT}^{\widering{T+U}}f[\vp_1(t)]\left|\zf\right|^2{\rm d}t= \\
& =\omega(T)\ln T\int_T^{T+U}f(x){\rm d}x,\ U\in\left(\left. 0,\frac{T}{\ln T}\right.\right],
\end{split}
\ee
where
\bdis
\omega(T)=\omega(T;f)=1+\mcal{O}\left(\frac{\ln\ln T}{\ln T}\right),\ T\to\infty.
\edis
Next, from (a) -- (c) of Theorem we obtain that
\be \label{7.2}
\begin{split}
& \int_{\gamma'}^{\gamma''} |G(t)|^\alpha |G'(t)|{\rm d}t=\frac{2}{\alpha+1}|G(t_0)|^{\alpha+1},\ \alpha>0, \\
& G(t)\geq 0 \ (\leq 0),\ t\in [\gamma',\gamma''].
\end{split}
\ee
Now, in the case
\bdis
f(t)=|G(t)|^\alpha |G'(t)|,\ t\in [\gamma',\gamma'']
\edis
we have (see (\ref{7.1}), (\ref{7.2}) and (b) of Theorem)
\be \label{7.3}
\begin{split}
& \int_{\mathring{\gamma}'}^{\mathring{\gamma}''}
|G[\vp_1(t)]|^\alpha |G'_{\vp_1}[\vp_1(t)]|\left|\zf\right|^2{\rm d}t= \\
& = \frac{2\omega(\gamma')\ln\gamma'}{\alpha+1}|G(t_0)|^{\alpha+1}.
\end{split}
\ee
Hence, putting
\bdis
\alpha=2\omega(\gamma')\ln\gamma'-1,
\edis
we obtain by (\ref{7.3}) and (e) of Theorem that
\be \label{7.4}
\begin{split}
&\int_{\mathring{\gamma}'}^{\mathring{\gamma}''}
|H[\vp_1(t);t_0]|^{2\omega(\gamma')\ln\gamma'-1}|H'_{\vp_1}[\vp_1(t);t_0]|\left|\zf\right|^2{\rm d}t=1.
\end{split}
\ee
Finally, if we use the mean-value theorem in (\ref{7.4}), we obtain
\bdis
\begin{split}
& \left|\zeta\left(\frac 12+it_H\right)\right|^2\int_{\mathring{\gamma}'}^{\mathring{\gamma}''}
|H[\vp_1(t);t_0]|^{2\omega(\gamma')\ln\gamma'-1}|H'_{\vp_1}[\vp_1(t);t_0]|{\rm d}t=1,
\end{split}
\edis
i. e. the formula (\ref{2.2}) is verified.

\section{Concluding remarks: the Dirac property of the topological deformation $H_T$ at the point $\mathring{t}_0$}

Since
\bdis
\vp_1(\mathring{t}_0)=t_0 \ \Rightarrow \ H_T[\vp_1(\mathring{t}_0);t_0]=1
\edis
(see (\ref{6.1})), then we have (see (\ref{7.4}))
\be \label{8.1}
\begin{split}
& \int_{\mathring{\gamma}'}^{\mathring{\gamma}''}
|H[\vp_1(t);t_0]|^{2\omega(\gamma')\ln\gamma'-1}|H'_{\vp_1}[\vp_1(t);t_0]|\left|\zf\right|^2{\rm d}t= \\
& = |H_T[\vp_1(\mathring{t}_0);t_0]|^{2\omega(\gamma')\ln\gamma'-1}.
\end{split}
\ee
Hence, by using the Dirac $\delta$-function (see \cite{1}, pp. 58-61) we obtain
\be \label{8.2}
\begin{split}
&\int_{\mathring{\gamma}'}^{\mathring{\gamma}''}
|H[\vp_1(t);t_0]|^{2\omega(\gamma')\ln\gamma'-1}\delta(t-\mathring{t}_0){\rm d}t= \\
& =|H_T[\vp_1(\mathring{t}_0);t_0]|^{2\omega(\gamma')\ln\gamma'-1}.
\end{split}
\ee
The set
\bdis
(-\infty,\mathring{\gamma}')\cup (\mathring{\gamma}'',\infty)
\edis
is irrelevant, comp. \cite{1}, p. 59.

\begin{remark}
By (8.1), (8.2) we see that the continuous function (\emph{proper} function in terminology of Dirac)
\bdis
|H'_{\vp_1}[\vp_1(t);t_0]|\left|\zf\right|^2
\edis
acts at the point $\mathring{t}_0$ as the Dirac $\delta$-function (\emph{improper} function in terminology of Dirac, see \cite{1}, p. 58).
\end{remark}

\begin{remark}
It is sufficient to quote Dirac as the discoverer of the generalized functions.
\end{remark}

\thanks{I would like to thank Michal Demetrian for his help with electronic version of this paper.}


\begin{thebibliography}{29}
%
\bibitem{1}
P.A.M. Dirac, `\emph{The principles of quantum mechanics}`, Oxford, Clarendon Press, 1958.
%
\bibitem{2}
A.A. Karatsuba, `\emph{Complex analysis in number theory}`, CRC Press, Boca Raton, Ann Arbor, London, Tokyo, 1995.
%
\bibitem{3}
J.E. Littlewood, `Two notes on the Riemann zeta-function`, Proc. Cambr. Phil. Soc., 22 (1924), 234-242.
%
\bibitem{4}
J. Moser, `Some properties of the Riemann zeta-function on the critical line`, Acta Arith., 26 (1974), 33-39, (in Russian); arXiv: 0710.0943.
%
\bibitem{5}
J. Moser, `Jacob's ladders and the almost exact asymptotic representation of the Hardy-Littlewood integral', Math. Notes 88, 414-422 (2010),
arXiv: 0901.3937.
%
\bibitem{6}
J. Moser, `Jacob's ladders, the structure of the Hardy-Littlewood integral and some new class of nonlinear integral equations`,
Proc. Stek. Inst. 276, 208-221 (2011), arXiv: 1103.0359.
%
\bibitem{7}
J. Moser, `Riemann hypothesis and some infinite set of microscopic universes of the Einstein's type in the early period of the
evolution of the Universe`, arXiv: 1307.1095.
%
\bibitem{8}
E.C. Titchmarsh, `\emph{The theory of the Riemann zeta-function}` Clarendon Press, Oxford, 1951.
%
\bibitem{9}
F.G. Tricomi, `\emph{Differential equations}`, Blackie and Son Limited, 1961 (in Russian).
%
\bibitem{10}
F.G. Tricomi, `\emph{Lezioni sule equationi a derivate parziali}`, Gheroni, Torino, 1954 (in Russian).
%
\bibitem{11}
E.T. Whittaker, G.N. Watson, `\emph{A course of modern analysis}`, Cambridge, Cambridge Univ. Press, 1927.

\end{thebibliography}
\end{document}